\documentclass[imslayout,noinfoline,a4paper,12pt,twoside]{imsart}
\RequirePackage[OT1]{fontenc} 
\RequirePackage{amsthm,amsmath,amssymb,natbib} 

\startlocaldefs
\theoremstyle{plain}
\newtheorem{thm}{Theorem}[section]

\newtheorem{lem}{Lemma}[section]
\endlocaldefs

\renewcommand{\[}{\begin{eqnarray*}}
\renewcommand{\]}{\end{eqnarray*}}
\newcommand{\la}{\begin{eqnarray}}
\newcommand{\al}{\end{eqnarray}}

\renewcommand{\epsilon}{\varepsilon}
\renewcommand{\phi}{\varphi}

\newcommand{\N}{{\mathbb N}}
\newcommand{\R}{{\mathbb R}}

\newcommand{\cP}{{\mathcal P}}

\newcommand{\Prob}{\mathrm{Prob}}
\newcommand{\E}{{\mathbb E}}  
\newcommand{\CE}[2]{\E\pmb{\left(\right.} #1 \pmb{\left|\right.} 
#2 \pmb{\left.\right)}}

\renewcommand{\d}{{\,\text{\rm d}}}

\begin{document}

\begin{frontmatter}


\title{Optimal L$_{}^1$-bounds for submartingales}

\runtitle{Optimal L$_{}^1$-bounds for submartingales }

\begin{aug}
\author{\fnms{Lutz} \snm{Mattner}
\ead[label=e1]{mattner@uni-trier.de}}
\and
\author{\fnms{Uwe} \snm{R\"osler}
\ead[label=e2]{roesler@math.uni-kiel.de}}%

\runauthor{Lutz Mattner and Uwe R\"osler}

\affiliation{Universit\"at Trier and Christian-Albrechts-Universit\"at
  zu Kiel} 
 {\rm \today}\\
 {\footnotesize\tt \jobname.tex}

\address{Universit\"at Trier\\
Fachbereich  IV - Mathematik\\ 
54286 Trier \\
Germany\\
\printead{e1}}

\address{Christian-Albrechts-Universit\"at zu Kiel \\
     Mathematisches Seminar \\
     Ludewig-Meyn-Str. 4 \\
     D-24098 Kiel  \\
     Germany \\
 \printead{e2}}
\end{aug}

\begin{abstract}
 The optimal function  $f$ satisfying 
 $$
   \E |\sum_{1}^n X_i | \,\,\ge\,\,f(\E|X_1|,\ldots,\E|X_n|) 
 $$ 
 for every martingale $(X_1,X_1+X_2, \ldots,\sum_{i=1}^n X_i)$ is shown to be 
 given by
 $$
  f(a)  \,\,=\,\, \max\, \Big\{ 
  a_k-\sum_{i=1}^{k-1} a_i\Big\}_{k=1}^n  \cup
  \Big\{\frac {a_k}2\Big\} _{k=3}^n 
 $$ 
 for $a\in{[0,\infty[}^n_{}$.
 A similar result is obtained for submartingales 
 $(0,X_1,X_1+X_2,\ldots, \sum_{i=1}^n X_i)$.

 The optimality proofs use 
 a convex-analytic comparison lemma of independent interest.
\end{abstract}

\begin{keyword}[class=AMS]
\kwd[Primary ]{60G42, 60E15}
\kwd[; secondary ]{26B25}.
\end{keyword}
\begin{keyword}
 \kwd{Convexity} 
 \kwd{extremal problem}
 \kwd{martingale}
 \kwd{moment inequalities}.
\end{keyword}

\end{frontmatter}



\section{Introduction and main results} 
This article provides in two interesting special cases the optimal
lower bounds for absolute moments of sums $S_n=\sum_{i=1}^n X_i$ given 
absolute moments of their increments $X_i$
and given some structural assumption on the process
$(X_1,\ldots,X_n)$, see Theorems~\ref{Thm MG case}
and~\ref{Thm SMG case} below. These belong to the first few nontrivial results 
of their kind, 
despite a rather large literature on moment bounds in general.
Let us introduce some notation 
for  stating our results 
and for comparing them  with previous ones.

For $n\in\N:=\{1,2,3,\ldots\}$ we consider processes $(X_1,\ldots,X_n)$
of real-valued random variables and always put $S_k:= \sum_{i=1}^k X_i$
for $k\in \{0,\ldots,n\}$, so that in particular $S_0 =0$.   
We denote various structural assumptions on $(X_1,\ldots,X_n)$
by acronyms as follows, 
\[
 \text{IC:}&&\hspace{-2ex} \text{$X_1,\ldots,X_n$ 
  are independent and centred}\\ 
 \text{MG:}&&\hspace{-2ex} \text{$(S_1,\ldots, S_n) $ is a martingale} \\
\text{SMG:} &&\hspace{-2ex}\text{$(S_0,\ldots, S_n) $ is a 
 submartingale  or a supermartingale}
\]
where  ``centred'' means $\E X_i = 0$  for each $i$.

For $r\in[1,\infty[$, 
$n\in\N$, 
$a\in{[0,\infty[}^n_{}$, and with  A denoting any of the assumptions 
IC, MG or SMG,  we put 
\la          \label{Def fA}
  f^{}_{r,\mathrm{A}}(a)
 &:=& \inf \,\{ \E |S_n |^r_{} \,:\,
     \text{\rm A},\, \E|X_1|^r_{} = a_1,\ldots,\E|X_n|^r_{}=a_n  \} \\ 
\label{Def FA}
  F^{}_{r,\mathrm{A}}(a)
 &:=& \sup \,\{ \E |S_n |^r_{} \,:\,
     \text{\rm A},\, \E|X_1|^r_{} = a_1,\ldots,\E|X_n|^r_{}=a_n  \}
\al
We sometimes abbreviate $f^{}_{1,\mathrm{A}}=:f^{}_\mathrm{A}$.

Thus, for example, $f^{}_\mathrm{MG}(a)$  with $a\in[0,\infty[^n$
is by definition the best lower bound for $\E |S_n|$ given 
that $(S_1,\ldots,S_n)$ is a martingale with  
$X_i=S_i-S_{i-1}$ satisfying $\E |X_i|=a_i$
for each $i\in\{1,\ldots,n\}$.

Set theoretical correctness requires us to note that   
the above definitions~\eqref{Def fA} and \eqref{Def FA},
referring as they do    to the nonexistent set of all processes
satisfying assumption~$\mathrm{A}$, can be legalized by 
noting that the expectations and assumptions considered 
depend only on the laws  of the processes $(X_1,\ldots,X_n)$.
These laws form a well-defined set on which a formally correct version
of definitions~\eqref{Def fA} and \eqref{Def FA} should  be based,
as in formula \eqref{ReDef.fn} below, which correctly defines 
the restriction of $f^{}_\mathrm{MG}$ to ${[0,\infty[}^n$.

The explicit  determination of the functions $f^{}_\mathrm{MG}$
and $f^{}_\mathrm{SMG}$ in the following two theorems is
the main  result of this paper:

\begin{thm}[The martingale case] \label{Thm MG case}
For $n\in \N$ and $a\in{[0,\infty[}^n_{}$, we have
\la         \label{Formel.f.MG}
f^{}_\mathrm{MG}(a)&=&
\max\, \Big\{   a_k-\sum_{i=1}^{k-1} a_i \Big\}_{k=1}^n \cup
\Big\{ \frac {a_k}2\Big\}_{k=3}^n 
\al
\end{thm}

\begin{thm}[The submartingale case] \label{Thm SMG case}
For $n\in \N$ and $a\in{[0,\infty[}^n$, we have 
\la
&&f^{}_\mathrm{SMG}(a)\,\,=\,  \label{Formel.f.SMG}\\ \nonumber
&&\phantom{,}\max\,\Big\{ a_k-\sum_{i\neq k} a_i \Big\}_{k=1}^n \cup
           \Big\{    
               \frac {a_k-\sum\limits_{i< k} a_i}2 \Big\}_{k=1}^{n-1} \cup 
               \Big\{      
                    \frac{a_k-\sum\limits_{i>k} a_i}3 \Big\}_{k=3}^{n} \cup 
              \Big\{    \frac {a_k}4 \Big\}_{k=3}^{n-1} 
\al
\end{thm}
Here and below $\{b_k\}_{k=n_1}^{n_2}$ denotes the possibly empty set
$\{b_k \,:\, k\in\N, \, n_1\le k \le n_2\}$ and summations as in 
\eqref{Formel.f.SMG} are with respect to 
$i\in\{1,\ldots,n\}$ subject to the indicated conditions.  
Put  less formally, Theorem~\ref{Thm MG case}, for example,  says
that inequality~\eqref{Expl.MG.bound}
below is optimal if just MG is assumed and if the right hand side
is only allowed to depend on $\E|X_1|,\ldots,\E |X_n|$.

Theorems~\ref{Thm MG case} and \ref{Thm SMG case} are proved in 
Section~\ref{Proofs} 
using in particular a comparison lemma for certain convex functions,
namely Lemma~\ref{Lem.convex.comparison}
supplied in the preparatory Section~\ref{Aux.Conv}, 
which might be of independent interest.
The remainder 
of this section contains a few  remarks
and some comparisons with known results.

In Theorem~\ref{Thm MG case} we could have written 
$\{ \frac {a_k}2\}_{k=1}^n$  instead of $\{ \frac {a_k}2\}_{k=3}^n$
either by changing the proof a bit  or by noting 
that $\frac {a_1}2 \le a_1 -\sum_{i=1}^{1-1} a_i$
and $\frac {a_2}2 = \frac 12 (a_1+a_2-a_1)\le \max\{a_1,a_2-a_1 \}$.
However, in the formulation chosen, the restriction of  $f^{}_\mathrm{MG}$
to ${[0,\infty[}^n$ is presented as a pointwise supremum of 
$n+(n-2)_+$ linear functions of which no proper subset would suffice.
(To check this, one may consider $a=(a_i)_{i=1}^n$
with $a_i = \delta_{k,i}$ to prove that  $a_k- \sum_{i=1}^{k-1}a_i$
can not be omitted for $k\in\{1,\ldots,n\}$,
and $a_i = \delta_{1,i} + \delta_{2,i} + 3 \delta_{k,i}$ 
to show  that  $\frac {a_k}2$
can not be omitted for $k\in\{3,\ldots,n\}$.)
A similar remark applies to  Theorem~\ref{Thm SMG case}.
(For proving here that none of the linear functions 
can be omitted,  one  may consider 
$a_i = \delta_{k,i}$ for $k\in\{1,\ldots,n\}$,
$a_i = 3\delta_{k,i} +  2 \delta_{n,i}$
for $k\in\{1,\ldots,n-1\}$,
$a_i = \delta_{1,i} + \delta_{2,i} + 2\delta_{k,i}$
for $k\in\{3,\ldots,n\}$,
and $a_i = \delta_{1,i} + \delta_{2,i} + 3\delta_{k,i}+ \delta_{n,i} $
for $k\in\{3,\ldots,n-1\}$.)
In particular, Theorems~\ref{Thm MG case} and \ref{Thm SMG case} 
 imply the inequalities
\la     \label{fMG.max}
 \frac 12  \max_{k=1}^n a_k &\le&  f^{}_\mathrm{MG}(a) 
    \,\,\le \,\,  \max_{k=1}^n a_k \label{fMG.max}     \\
 \frac 14  \max_{k=1}^n a_k &\le&  f^{}_\mathrm{SMG}(a)
 \,\,\le \,\,  \max_{k=1}^n a_k \label{fSMG.max}
\al
for $n\in\N$ and $a\in{[0,\infty[}^n$, where the constant factors
involved are easily seen to be optimal. See \eqref{lower bound gn} 
below for more details on \eqref{fSMG.max}.
Another consequence of  Theorems~\ref{Thm MG case}
is the Kemperman \& Smit~(1974, statement (d)) inequality 
\la   \label{K.S.Ineq}
  \frac 1{2n-1} \sum_{i=1}^n a_i &\le&  f^{}_\mathrm{MG}(a)
\al
since \eqref{Formel.f.MG} with 
$\{ \frac {a_k}2\}_{k=1}^n$  instead of $\{ \frac {a_k}2\}_{k=3}^n$
yields
$
 f^{}_\mathrm{MG}(a)
\ge \max \{a_1,\frac{a_2}2,\ldots,\frac{a_n}2 \} \linebreak
\ge  \frac 1{2n-1}a_1 + \sum_{k=2}^n \frac 2{2n-1}\, \frac{a_k}2
=\text{L.H.S.\eqref{K.S.Ineq}}$, where in the second step we 
have bounded a maximum from below by a convex combination.
However, depending on $a\in{[0,\infty[}^n$,
inequality~\eqref{K.S.Ineq} can be much worse 
but can only be slightly better than the left hand inequality 
in~\eqref{fMG.max}.

Theorem~\ref{Thm MG case} and its proof remain valid if 
we replace $(S_i)_{i=1}^n$ by  $(S_i)_{i=0}^n$
in the definition of $\mathrm{MG}$ and in the line 
following \eqref{The martingale}, and $\{2,\ldots,n\}$ by $\{1,\ldots,n\}$
in the line following \eqref{Canon.Cond.MG}.
No such remark applies to Theorem~\ref{Thm SMG case},
as becomes clear by considering the submartingale $(S_1,S_2):=(-1,0)$.

We now present all other cases known to us where one of the functions
$ f^{}_{r,\mathrm{A}}$ or $F^{}_{r,\mathrm{A}}$ with $r\in{[1,\infty[}$ and 
$\mathrm{A} \in\{\mathrm{IC}, \mathrm{MG}, \mathrm{SMG}\}$ is obvious 
or has been determined in the literature. For comparison we also mention 
some corresponding results referring to one of the assumptions 
\[
 \text{IS:}&&\hspace{-2ex} \text{$X_1,\ldots,X_n$ 
  are independent and symmetrically distributed}\\ 
 \text{IIDC:}&&\hspace{-2ex} \text{$X_1,\ldots,X_n$ 
  are independent, identically distributed, and centred} \\
\text{N:} &&\hspace{-2ex}\text{No assumption, i.e.,
$X_1,\ldots,X_n$ are arbitrary random variables}
\]
We make no attempt to review moment 
inequalities optimal in senses weaker than ours, like bounds of optimal 
order or optimal constants in bounds of a special form,
let alone related bounds involving moments of different orders or
tail probabilities. 

As a preparation, let $r\in{[1,\infty[}$ and let us consider 
$n\in\N$, $a\in{[0,\infty[}^n$,  $p\in{]0,1]}$,
and independent 
\[
 X_i&\sim&  (1-p)\delta_0 +\frac p2
 \big(\delta_{-a_i^{1/r}/p}+ \delta_{a_i^{1/r}/p}\big)
\]
for $i\in \{1,\ldots,n\}$.
Then $\E |X_i|^r= a_i$ for each $i$ and,
by bounding the probability of the event $\{X_i X_j\neq 0$ for some $i\neq j\}$,
we have $\lim_{p\downarrow 0}  \E |S_n|^r 
= \lim_{p\downarrow 0}  \E\sum_{i=1}^n |X_i|^r
= \sum_{i=1}^n a_i =: \|a\|_1$. 
This yields
\la \label{f,ell1,F}
 f^{}_{r,\mathrm{A}} & \le&   \|\cdot \|_1  \,\,\le \,\, F^{}_{r,\mathrm{A}}
\al
whenever $\mathrm{A}$ is any of the six  assumptions introduced above,
if in case of $\mathrm{A}= \text{IIDC}$ attention is restricted 
to arguments $a\in{[0,\infty[}^n$  with $a_1 = \ldots =a_n$.

Continuing now with the exponent $r=1$ but turning to upper bounds,
we  have 
$F^{}_{1,\mathrm{IC}} = F^{}_{1,\mathrm{MG}}=F^{}_{1,\mathrm{SMG}}=\|\cdot\|_1$
by the $\mathrm{L}^1$ triangle inequality 
$\E |S_n| \le \sum_{i=1}^n \E |X_i|$ 
and by the right hand inequality in \eqref{f,ell1,F} with $r=1$.

Coming back to lower bounds for $r=1$, we believe that $f^{}_{1,\mathrm{IC}}$
is in general unknown, although it is an easy exercise to show 
that   $f^{}_{1,\mathrm{IC}}(a) = \max_{k=1}^n a_k$ whenever 
$n\in \N$ and $a\in{[0,\infty[^n}$ satisfy $a_k \ge \sum_{i\neq k} a_i$ 
for some $k$. 
Somewhat surprisingly,  the corresponding problem under the assumption 
$\mathrm{IIDC}$
has  already been solved in  Mattner~(2003), 
namely  $f^{}_{1,\mathrm{IIDC}}(a)= c_n a_1$ for $a\in{[0,\infty[}^n$ with 
$a_1=\ldots=a_n$, where $c_n$ is $n$ times the maximal probability of a 
binomial distribution with parameters $n$ and $\lfloor n/2\rfloor/n$,
and thus  $ c_n\sim \sqrt{\frac {2n}{\pi}}$ for 
$n\rightarrow \infty$.

For  $r=2$, we obviously have  
$f^{}_{2,\mathrm{IC}}=F^{}_{2,\mathrm{IC}}=f^{}_{2,\mathrm{MG}}=F^{}_{2,\mathrm{MG}}
= \|\cdot\|_1$.

Let $r\in {]1,2[}$. Then  von Bahr \& Esseen~(1965, Theorem 1) 
and the right hand inequality in~\eqref{f,ell1,F} yield
$F_{r,\mathrm{IS}} = \|\cdot\|_1$. 

Let now $r\in [3,\infty[$. Then $\E |X+Y|^r \ge \E |X|^r + \E |Y|^r$
whenever  $X,Y$ are independent and centred $\R$-valued random variables,
by Cox \& Kemperman~(1983, Theorem 2.6). This immediately yields
the remarkable result
 $f^{}_{r,\mathrm{IC}} = \|\cdot\|_1 $,
by induction and by the left hand inequality in 
 \eqref{f,ell1,F}.
See  Pinelis~(2002, Theorem 6) for an extension of the Cox-Kemperman 
inequality  to Hilbert space valued random variables.

Finally, computing   $f_{r,\mathrm{N}}$ and $F_{r,\mathrm{N}}$
for $r\in[1,\infty[$  is equivalent to determining the best lower 
and upper bounds for the $\mathrm{L}^r$ norm
of a sum given the norms of the summands. This is a special case
of an exercise in elementary normed vector space theory solved 
in  Mattner~(2008). The solution given there yields
\[
 f_{r,\mathrm{N}}(a) &=&\max_{k=1}^n  
 \Big( (a_k)^{\frac 1r}_{} - \sum_{i\neq k} (a_i)^{\frac 1r}_{}\Big)^{r}_{+}\\
 F_{r,\mathrm{N}}(a) &=&   \Big(\sum_{i=1}^n (a_i)^{\frac 1r}_{} \Big)^{r}_{} 
\]
for $n\in\N$ and $a\in{[0,\infty[}^n$.

\section{Auxiliary facts from convex analysis}
\label{Aux.Conv}
We assume as known  some standard terminology and facts as 
given in Rockafellar~(1970), but generalized in the obvious way from $\R^n$
to arbitrary and possibly infinite-di\-men\-sio\-nal
vector spaces over $\R$. 
We state here two lemmas used in
Section~\ref{Proofs} below.

The first lemma says that partial infima of convex functions are
convex.  
This  will be applied below only in situations where 
all functions considered 
are finite-valued, but it seems simpler to state the general case.  
\begin{lem}\label{Lem.partial.inf.convex}
Let $C$ and $D$ be convex subsets of vector spaces over $\R$ and let
$E\subset C\times D$ and $f:E\rightarrow [-\infty,\infty]$ be convex.
Then the function
\la      \label{partial.inf} \nonumber
 C\,\ni\, x&\mapsto& \inf\{ f(x,y)\,:\,(x,y)\in E \} \,\in\,[-\infty,\infty]  
\al
is convex.
\end{lem}
\begin{proof}
Easy and well-known, compare Rockafellar~(1970, pp.~38-39). 
\end{proof} 

Our second lemma reduces the pointwise comparison of certain convex functions
to a comparison at  distinguished  points. Let $E$ be a vector space over $\R$.
Then the  {\em dimension} of a convex set  $C\subset E$ is, 
by definition,
the dimension of the vector subspace $F\subset E$ obtained by translating the
affine hull $A$ of $C$ towards the origin.  
If this dimension is finite,
then $A$ is topologized  by  translating from $F$ its 
usual (unique Hausdorff topological vector space) topology,
and the {\em relative boundary} and the {\em relative interior}
of $C$ are  then the boundary and the interior 
of $C$ as a subset of the topological space $A$.

\begin{lem}\label{Lem.convex.comparison} 
Let $C\subset E$ be  convex, finite-dimensional, and  compact.
Let $f,g:C \rightarrow \R$ be functions with $f$ convex and 
$g=\sup_{i\in I}g_i $ being the pointwise supremum of a finite family
of affine functions $g_i: C \rightarrow \R$. Assume that $f(x)\le g(x)$
holds for every $x$ satisfying one   of the two conditions
\la
&&x \,\in \text{ \rm  relative boundary of } C      \label{x.rel.bound.C}  \\
&&x \,\in \text{ \rm  relative interior of } C \,\,\text{ \rm and  }\,\,
     \{g_i\,:\, i\in I,\,g_i(x)=g(x) \}\text{ \rm contains} \label{x...k+1} \\
  && \nonumber \text{\rm  at least
       $\dim C + 1$ affinely independent functions}  
\al
Then $f\le g$ on $C$. 
\end{lem}

We recall that a family $(g_i:i\in J)$ is {\em affinely independent}
if $\sum_{i\in J}^n \alpha_i g_i=0$ with $\alpha_i\in\R$ and 
$\sum_{i\in J } \alpha_i=0$ implies $\alpha_i=0$ for $i\in J$.

For example, to prove $f(x):= x^2 \le |x| =: g(x)$ for $x\in [-1,1]=:C$,
it is enough to consider $x\in\{-1,1,0\}$, by 
Lemma~\ref{Lem.convex.comparison} 
applied with $g_1(x):=-x$ and  $g_2(x):=x$. Replacing here  $x^2$ by
$\frac 12 (x^2 +1)$ shows that the affine independence 
requirement in~\eqref{x...k+1} can not be strengthened to linear
independence. 

\begin{proof}[Proof of Lemma \ref{Lem.convex.comparison}]
We may exclude the trivial cases where $I$ is empty or $C$ has at most
one element. Thus $k:=\dim C \ge 1$.   
By a translation and by choosing a basis,
we may further assume that $C\subset \R^k=E$.
Then $C$ has nonempty interior $\mathrm{int\,} C$ as a subset of $\R^k$, 
and we can omit the adjective
``relative'' in conditions \eqref{x.rel.bound.C} and 
\eqref{x...k+1}. 
Also, we can consider the affine functions $g_i$
as being defined on the entire space~$\R^k$.

For every $x\in C$,  we introduce the nonempty set  
\[
 I_x  &:=&\{ i\in I: g_i(x)=g(x)\}
\]

Let $x_0 \in C$. Then 
\[
 x_0 \,\,\in\,\,  C_0 
 &:=& \{ x\in C: g_i(x)= g(x) \text{ for } i\in I_{x^{}_0} \}
\]
So it suffices to prove  $f\le g$ on $C_0$.  
Now $C_0$ is easily seen to be 
convex and compact, $f$ is convex, and 
$g$ is affine on $C_0$ as $I_{x_0}\neq \emptyset$. Hence,
by 
Rockafellar~(1970, Corollary 18.5.1),
it suffices to prove $f(x) \le g(x)$ whenever $x$ is an extreme point
of $C_0$. We finish this proof by showing that 
every extreme point $x$ of $C_0$ satisfies \eqref{x.rel.bound.C} or
\eqref{x...k+1}:

Let $x\in C_0$ satisfy neither \eqref{x.rel.bound.C} nor
\eqref{x...k+1}. Then $x\in \mathrm{int\,}C$
and $\{g_i: i\in I_x\}$ contains 
at most $k$ affinely independent functions. 
Also, $g_i(x)<g(x)$ for 
$i\in I\setminus I_x$  and thus 
\[
x\,\,\in\,\,  U &:=& 
  \big(\text{\rm int\,} C \big)\cap \{y\in\R^k: g_i(y)<g(y) 
    \text{ for } i\in I\setminus I_x \}
\]
and $U$ is open, due to the finiteness of $I$.
Let us assume for notational convenience that 
$I_x = \{1,\ldots,n\}$.  Assuming $n\ge 2$ until further notice, we have 
\[
 x\,\,\in\,\, 
  A&:=& \{y\in\R^k: g_1(y)=\ldots=g_n(y)\} \\
  &=&  \{y\in\R^k: h_1(y)=0, \ldots, h_{n-1}(y)=0\} \\
  &=& \{y\in\R^k: h_1(y)=0, \ldots, h_{\ell}(y)=0\}
\]
where $h_i := g_i-g_n$ for $i\in\{1,\ldots,n-1\}$ and,
reordering if necessary,  
$(h_1,\ldots,h_\ell)$ is a maximal linearly independent 
subfamily of $ (h_1,\ldots,h_{n-1})$. Then 
$(g_1,\ldots,g_\ell, g_n)$ is affinely independent,
so $\ell  +1 \le k $.
As the functions $h_i$ are affine and the set $A$ is nonempty,
it follows that  
$\dim A \ge k- \ell \ge 1$. 
Thus $U\cap A$ contains a nondegenerate line segment $S$ 
with midpoint $x$, this conclusion is also true if $n=1$ and $A:= \R^k$,
and we no longer assume $n\ge 2$. 
For $y\in U\cap A$,
we have $g_i(y) =g(y) $ for every $i\in I_x$.
As $I_x \supset I_{x_0}$, we conclude that $ U\cap A \subset C_0$
and hence $S\subset C_0$, so that $x$ is not an extreme point
of~$C_0$. 
\end{proof}
\section{Proofs of the main results} \label{Proofs}
\begin{proof}[Proof of Theorem \ref{Thm MG case}] 
In this proof we put 
\[
 f(a) \,\,:=\,\,f^{}_\mathrm{MG}(a)  \,\,=\,\,\text{L.H.S.\eqref{Formel.f.MG}} 
 &\text{and}& g (a) \,\,:=\,\, \text{R.H.S.\eqref{Formel.f.MG}}
\]
for $n\in \N$ and $a\in[0,\infty[^n$. 

Let $n\in \N$ and let $(S_k)_{k=1}^n = (\sum_{i=1}^k X_i)_{k=1}^n$
be a martingale. Then
\la \label{Ineq.ESk.MG}
  \E | S_k| &\le& \E|S_n| \qquad\text{for }k\in\{1,\ldots,n\}
\al
Applying first the $\mathrm{L}^1$ triangle inequality and
then \eqref{Ineq.ESk.MG} to each of the identities
\[\begin{array}{rcll}
\displaystyle X_k &=&\displaystyle S_k -  \sum_{i=1}^{k-1} X_i 
 &\quad\text{for }k\in\{1,\ldots,n\}\\
\displaystyle X_k &=&\displaystyle S_k - S_{k-1}\phantom{\sum^{}} 
 &\quad\text{for }k\in\{3,\ldots,n\}
\end{array}\]
yields $ \E |X_k| \le \E |S_n| + \sum_{i=1}^{k-1} \E |X_i|$
and  $ \E |X_k| \le 2 \E |S_n|$, respectively, and hence
\la      \label{Expl.MG.bound}
   \E |S_n| &\ge& 
 \max\, \Big\{\E |X_k| -\sum_{i=1}^{k-1} \E |X_i|\Big\}_{k=1}^n \cup 
\Big\{ \frac {\E |X_k|}2\Big\}_{k=3}^n 
\al
This proves $f\ge g$.

It remains to prove the reversed inequality  $f\le g$, and this will 
eventually be done  by induction.  For $n\in \N$, we let $f_n$ and $g_n$ denote 
the restrictions of $f$ and $g$ to ${[0,\infty[}^n$.   

A key observation is that each $f_n$ is convex. To see this,
let us consider the canonical process 
$X=(X_1,\ldots,X_n):=\mathrm{id}_{\R^n}$ and the set of laws
\[
\cP &:=& \{P \in\Prob(\R^n) : X\text{ satisfies MG under }P\}
\]
The latter is convex, since 
$P\in \Prob(\R^n)$ belongs to $\cP$ iff 
it satisfies the linear constraint
\la     \label{Canon.Cond.MG}
  \int_{\R^n} X_i\, h(X_1,\ldots,X_{i-1})\d P &=&0
\al
for each $i\in \{2,\ldots,n\}$ and each measurable
indicator $h:\R^{i-1}\rightarrow \{0,1\}$. 
Recalling the notation $S_n:=\sum_{i=1}^n X_i$, we have 
\la     \label{ReDef.fn}\quad
  f_n(a)&=&\inf\Big\{ \int_{\R^n} \left|S_n \right|\d P\,:\, 
   P\in \cP,\, \int_{\R^n}\left| X_i\right|\d P  =a_i 
  \text{ for each }i
 \Big\}
\al
for $a\in{[0,\infty[}^n$.
Thus an application of Lemma~\ref{Lem.partial.inf.convex},
with  $C= {[0,\infty[}^n$,
$D=\cP$, 
$E =  \{ (a,P) \in C\times D \,:\,\int_{\R^n}\left| X_i\right|\d P  =a_i 
  \text{ for each }i \}$
and 
$ f = ( (a,P) \mapsto   \int_{\R^n} \left|S_n \right|\d P)$,
yields the claimed convexity of $f_n$. 

Next, the functions $f_n$ and $g_n$ are homogeneous: 
Since constant multiples of martingales are martingales, we have 
$f_n(\lambda a_1,\ldots, \lambda a_n) = \lambda f_n(a_1,\ldots, a_n)$
for $\lambda \in [0,\infty[$ and $a\in {[0,\infty[}^n$,  
and the same scaling relations obviously hold for the functions $g_n$.

By the homogeneity just observed, it suffices to prove
$f_n\le g_n$ on the simplex 
$C_n:=\{ a\in{[0,\infty[}^n : \sum_{i=1}^n a_i =1\}$, for each $n\in \N$.
The case $n=1$ is trivial. So let us assume that 
$n\in \N$ with $n\ge 2$ and that $f_{n-1} \le g_{n-1}$ on $C_{n-1}$.
To prove $f_n \le g_n$ on $C_n$, we will apply 
Lemma~\ref{Lem.convex.comparison}  with $C=C_n$,
$\dim C =n-1$, $f=f_n$, $g=g_n$, 
and $(g_i:i\in I) = (g^{\nu}_k \,:\, \nu\in\{1,2\},\,k\in K_\nu  ) $ 
with 
\[\begin{array}{rcll}
\displaystyle g^1_k(a)&:=&\displaystyle a_k-\sum_{i=1}^{k-1} a_i &\qquad
 \text{ for }k\in K_1:=\{1,\ldots,n\}\\  
\displaystyle g^2_k(a)&:=&\displaystyle \frac{a_k}2 \qquad\quad\;
 &\qquad\text{ for }k\in K_2:=\{3,\ldots,n\}
\end{array}\]
for $a\in{[0,\infty[}^n$. 

Suppose first that $a\in C_n$ belongs to the
relative boundary of $C_n$. Then for some  
$j \in \{1,\ldots,n\}$ we  have
$a_j=0$ and hence
\[
\begin{array}{rcll}
 f_n(a)&=& f_{n-1}(a_1,\ldots,a_{j -1},a_{j +1},\ldots,a_n)
             &\quad [\text{obvious by definition of $f_n$}] \\
 &\le& g_{n-1}(a_1,\ldots,a_{j -1},a_{j +1},\ldots,a_n)
             &\quad [\text{by induction hypothesis}]\\
 &=& g_n(a) & \quad  
 [\text{obvious by definition of $g_n$}] \\
\end{array}
\]

Now suppose that $a\in C_n$ is as in \eqref{x...k+1}, that is,
$a$ belongs to the
relative interior of $C_n$ and with 
\la \label{DefK.nu.1.2}
 K_\nu(a) &:=& \{ k \in K_\nu\,:\, g^\nu_k(a) =g_n(a)\}
  \qquad \text{ for }\nu\in\{1,2\}
\al
we have at least $\dim C_n +1=n$ affinely independent functions
in $\{g^{\nu}_k \,:\, \nu\in\{1,2\},\,k\in K_\nu(a)\}$, 
so that in particular
\la \label{card Ia}
  \sharp K_1(a) +  \sharp K_2(a) &\ge &  n \,\,\ge\,\, 2
\al
and $a_i >0$ for every $i\in\{1,\ldots,n\}$. We now must have 
\la  \label{K1(a),K2(a)}
  K_1(a)= \{1,2\},&& K_2(a) = \{3,\ldots,n\}
\al
for otherwise, in view of \eqref{card Ia} and $\sharp K_2 = n-2$,
there would exist $k,\ell \in K_1(a)$ with $k<\ell$ and $\ell \ge 3$,
so that 
\[
 a_\ell -\sum_{i=1}^{\ell -1} a_i 
 &=& g^1_\ell (a) \,=\, 
 g_n(a)  \,=\, 
 g^1_k(a) 
 \,=\, a_k -\sum_{i=1}^{k-1} a_i
\,=\, 2a_k -\sum_{i=1}^{k } a_i
\]  
yielding 
\[
 a_\ell &=& 2a_k + \sum_{i=k+1}^{\ell -1} a_i 
 \,\,>\,\, 2 \,\Big(a_k -\sum_{i=1}^{k-1} a_i  \Big)
\]
since one of the two sums  above is nonempty,
and hence the contradiction 
\[
  g_n(a) &=&  g^1_k(a)\,\,  = \,\,  a_k -\sum_{i=1}^{k -1} a_i
 \,\,  <\,\,   \frac{a_\ell}2\,\,  =\,\,  g^2_\ell(a) \,\,  \le \,\,  g_n(a)
\]

By \eqref{K1(a),K2(a)} we have 
$g_n(a)= a_1 = a_2-a_1 = \frac {a_3}2 =\ldots =\frac{a_n}2$
and thus 
\la   \label{a to be checked, MG}
 a_1 = \frac {a_2}2 =\ldots =\frac{a_n}2 
\al
To prove $f_n(a) \le g_n(a)$ in this case, let us consider
$p\in {]0,1[}$ and  
independent random variables $Y_1,\ldots,Y_n$ 
with 
\[
 Y_1 \sim \frac 12 \delta_{-1}+\frac 12 \delta_1,
 &&Y_i\sim (1-p)\delta_0 + p \delta_{1/p}
 \qquad \text{ for }i\in\{2,\ldots,n\}
\]
and let us put 
\la       \label{The martingale}
&&S_0 := 0, \quad S_k := \prod_{i=1}^k Y_i,\quad 
 X_k := S_k-S_{k-1} \quad\text{ for } k\in\{1,\ldots,n\}
\al
Then, since $\E Y_i=1$ for $i\ge 2$, the process
$(S_i)_{i=1}^n $ is a martingale.
We have $\E |S_i| =1$ for $i\in  \{1,\ldots,n\}$ and hence  $\E |X_1| =1$
and 
\[
 \E|X_i| &=& \E |Y_i-1|\E|S_{i-1}| 
 \,=\, \big((1-p)+p(\frac 1{p} -1) \big)\cdot 1 = 2-2p 
\]
for $i\in\{2,\ldots,n\}$. Thus 
\[
  f_n(1,2-2p,\ldots,2-2p) &\le & 1
\]
By the continuity of the convex function $f_n$  on the open set
${]0,\infty[}^n$, compare Rockafellar~(1970, p.~82), 
we can let $p\rightarrow 0$ to deduce $f_n(1,2,\ldots,2) \le 1$.
Hence,  by the homogeneity of $f_n$, we get
\[
 f_n(a) &\le & a_1 \,=\, g_n(a)
\]
for the $a$ satisfying \eqref{a to be checked, MG}. 

By the above, an application of Lemma~\ref{Lem.convex.comparison} 
yields $f_n\le g_n$ on $C_n$, which completes our 
inductive proof of $f\le g$.
\end{proof}

\begin{proof}[Proof of Theorem \ref{Thm SMG case}] 
This is parallel to but more complicated than the previous proof
we assume the reader has studied. We may and do replace 
``submartingale or supermartingale'' in the definition of $\mathrm{SMG}$ 
by ``submartingale''. We then put 
\[
 f(a) \,\,:=\,\,f^{}_\mathrm{SMG}(a)  \,\,=\,\,\text{L.H.S.\eqref{Formel.f.SMG}} 
 &\text{and}& g (a) \,\,:=\,\, \text{R.H.S.\eqref{Formel.f.SMG}}
\]
for $n\in \N$ and $a\in{[0,\infty[}^n$.

Let $n\in \N$ and let $(S_k)_{k=0}^n = (\sum_{i=1}^k X_i)_{k=0}^n$
be a submartingale. Then 
\la \label{Ineq.ESk.SMG}
  \E | S_k| &\le& 2 \E|S_n| \qquad\text{for }k\in\{1,\ldots,n-1\}
\al
by  Doob (1953, page 311, Theorem 3.1 (ii)), or by noting
that $(0,S_k,S_n)$ is a submartingale,  that is, 
$\E S_k\ge 0$ and $\CE{S_n}{S_k}\ge S_k$, and so 
$\E |S_k| \le \E |S_k| + \E S_k = 2 \E (S_k)^{}_+ \le 
2\E (\CE{S_n}{S_k})_+ 
\le 2\E \CE{(S_n)^{}_+}{S_k} \le 2 \E |S_n|$. 

Applying now first the $\mathrm{L}^1$ triangle inequality to each of the 
identities 
\[
\begin{array}{rcll}
 X_k &=&\displaystyle S_n- \sum_{i\neq k} X_i 
  &\qquad \text{for }k\in\{1,\ldots,n\}\\
 X_k &=&\displaystyle S_k- \sum_{i< k}^{} X_i 
  &\qquad\text{for } k\in\{1,\ldots,n-1\}\\
 X_k &=&\displaystyle S_n -S_{k-1} - \sum_{i> k}^{} X_i
  &\qquad\text{for }k\in\{3,\ldots,n\}\\
 X_k &=&\displaystyle S_k-S_{k-1}\phantom{\sum^{}}  &\qquad\text{for } k\in\{3,\ldots,n-1\}
\end{array}
\]
and then \eqref{Ineq.ESk.SMG} to the
resulting inequalities in  the last  three groups yields
$\E |X_k| \le \E|S_n| + \sum_{i\neq k} \E|X_i|$, 
$\E |X_k| \le 2\E|S_n| + \sum_{i< k} \E|X_i|$,
$\E |X_k| \le 3\E|S_n| + \sum_{i> k} \E|X_i|$,
and $\E |X_k| \le 4\E|S_n|$, respectively, and hence 
the analogue to \eqref{Expl.MG.bound} proving $f\ge g$ in the present case.

To prove $f\le g$, let $f_n$ and $g_n$ denote the restrictions
of $f$ and $g$ to ${[0,\infty[}^n$. 

We have
\la                  \label{lower bound gn} 
 g_n(a) &\ge& \max\, \big\{\frac{a_1}2\big\} 
\cup\big\{\frac{a_k}4\big\}_{k=2}^{n-1}
 \cup \big\{\frac  {a_n}3 \big\} 
 \,\,\ge\,\,\frac 14 \max_{k=1}^na_k
\al
for $n\in \N$ and $a\in[0,\infty[^n$, trivially if $n=1$ and otherwise
since $ g_n(a) \ge \frac{a_1}2$,  
$ g_n(a) \ge \max\{\frac{a_1}2,\frac{a_2-a_1}2\} \ge \frac{a_2}4$ if $n\ge 3$,
$ g_n(a) \ge \frac {a_k}4$ if $k\in\{3,\ldots,n-1\}$,
$ g_2(a) \ge \max\{a_2-a_1,\frac{a_1}2 \} \ge \frac{a_2}3$,
and $ g_n(a) \ge \frac{a_n}3$ if $n\ge 3$.

Each $f_n$ is convex
as in the previous proof, where we only have to replace
the equality sign in \eqref{Canon.Cond.MG}  by $\le$ and 
$\{2,\ldots,n\}$ by $\{1,\ldots,n\}$, and the functions 
$f_n$ and $g_n$ are homogeneous.
Proceeding again by induction, we apply
Lemma~\ref{Lem.convex.comparison} as above, with 
$(g_i:i\in I) = (g^{\nu}_k \,:\, \nu\in\{1,2,3,4\},\,k\in K_\nu  ) $ 
where
\[
\begin{array}{lcll}
 g^1_k(a)&:=&\displaystyle a_k-\sum_{i\neq k} a_i&\qquad
 \text{ for }k\in K_1:=\{1,\ldots,n\}\\ 
 g^2_k(a)&:=& \displaystyle\frac 12 \Big( a_k-\sum_{i<k} a_i \Big)& \qquad
 \text{ for }k\in K_2:=\{1,\ldots,n-1\}\\ 
 g^3_k(a)&:=& \displaystyle\frac 13 \Big( a_k-\sum_{i> k} a_i \Big)& \qquad
 \text{ for }k\in K_3:=\{3,\ldots,n\}\\ 
 g^4_k(a)&:=&\displaystyle \frac{a_k}4& \qquad
 \text{ for }k\in K_4:=\{3,\ldots,n-1\}
\end{array}
\]
for $a\in{[0,\infty[}^n$. Trivially, $f_1 \le g_1$.
So let $n\in \N$ with $n\ge 2$. As above, the induction hypothesis 
$g_{n-1}\le f_{n-1}$ yields $f_n(a) \le g_n(a) $ for $a$ belonging to 
the relative boundary of $C_n$.  Hence, defining $K_\nu(a)$ as 
in~\eqref{DefK.nu.1.2}  but now with 
$\nu\in\{1,2,3,4\}$, we assume for the rest of this proof that 
$a\in{]0,\infty[}^n$ and that 
we have at least $\dim C_n +1=n$ affinely independent functions
in $\{g^{\nu}_k \,:\, \nu\in\{1,2,3,4\},\,k\in K_\nu(a)\}$, 
so that in particular
\la\label{card.K1.4}
  \sharp K_1(a) +  \sharp K_2(a) +  \sharp K_3(a)  +  \sharp K_4(a)
   &\ge & 
 n \,\,\ge \,\,2 
\al

We will now prove in twelve steps that one of the two conditions 
\la                           \label{a_1/2...a_n/3}
 g_n(a)&=&  \frac{a_1}2\,\, =\,\, \frac{a_2}4 \,\,=  \,\ldots \,
  =\,\, \frac{a_{n-1}}4 \,\,=\,\, \frac{a_n}3  \\
 g_n(a)&=&
 \frac{a_1}2\,\, =\,\, \frac{a_2}4\,\, = \, \ldots \,
=\,\, \frac{a_{n-1}}4\,\, = \,\,a_n   \label{a_1/2...a_n} 
\al%
holds.  
(For $n=2$, condition  
\eqref{a_1/2...a_n/3} reads $g_n(a)= \frac{a_1} 2=\frac{a_2}3$.
Similarly for \eqref{a_1/2...a_n}.)

\textsc{Step 1:} $g_n(a)>0$.

\textsc{Proof:}  $g_n(a)\ge g^2_1(a)=\frac{a_1}2 >0$.

\textsc{Step 2:} $\sharp K_1(a)  \le 1$.

\textsc{Proof:}  
Otherwise there are $k,\ell \in\{1,\ldots,n\}$
with $k\neq \ell$ and $g^1_k(a)=g^1_\ell(a)=g_n(a)$. 
The first equation reads 
\[
 a_k-\sum_{i\neq k}a_i &=& a_\ell -\sum_{i\neq \ell}a_i 
\]
and yields  $a_k - a_\ell =a_\ell -a_k$, 
hence $a_k- a_\ell =0$, contradicting Step 1 through
\[
 g_n(a) \,=\, g^1_k(a) \,=\,a_k-\sum_{i\neq k}a_i \,\le\, a_k-a_\ell
 \,=\, 0
\]

\textsc{Step 3:} $\sharp K_3(a)  \le 1$.

\textsc{Proof:} Otherwise  there are $k,\ell \in\{3,\ldots,n\}$
with $k< \ell$ and $g^3_k(a)=g^3_\ell(a)=g_n(a)$. 
Hence $a_k-\sum_{i>k} a_i = a_\ell - \sum_{i>\ell} a_i$ and thus
\[            
  a_k &=& a_\ell +\sum_{i=k+1}^\ell a_i \,\ge \,\, 2a_\ell 
\]
yielding  the contradiction
\[
&\displaystyle
  g_n(a)\,=\, g^3_\ell(a) \,=\, \frac 13\Big(a_\ell -\sum_{i>\ell} a_i \Big)
 \,\le \, \frac {a_\ell}3 \,<\, \frac {a_\ell}2 \,\le \, \frac {a_k}4
 \,\,=\,\, g^4_k(a) 
 \,\le \, g_n(a) 
\] 

\textsc{Step 4:} $\sharp K_2(a)  \ge 1$.

\textsc{Proof:} Otherwise, by Steps 2 and 3, we would have
$\text{L.H.S.\eqref{card.K1.4}} \le 1+0+\min\{1,n-2\}+\max\{0,n-3\}=n-1$.

\textsc{Step 5:} \textsl{If $\sharp K_2(a)  \ge 2$, then $K_2(a)=\{1,2\}$}.

\textsc{Proof:} Otherwise there are $k,\ell \in\{1,\ldots,n-1\}$
with $k<\ell$, $\ell \ge 3$, and $g^2_k(a)=g^2_\ell(a)=g_n(a)$. 
Hence $a_\ell -\sum_{i<\ell}a_i = a_k -\sum_{i<k}a_i = 2a_k -\sum_{i\le k}a_i$ 
yielding 
\[
 a_\ell &=& 2a_k + \sum_{i=k+1}^{\ell -1} a_i 
 \,\,>\,\, 2 \,\Big(a_k -\sum_{i=1}^{k-1} a_i  \Big)
\]
since one of the two sums  above is nonempty,
and hence the contradiction 
\[
  g_n(a) &=&  g^2_k(a)\,\,  = \,\, 
  \frac 12 \Big( a_k -\sum_{i=1}^{k -1} a_i\Big)
 \,\,  <\,\,   \frac{a_\ell}4\,\,  =\,\,  g^4_\ell(a) \,\,  \le \,\,  g_n(a)
\]

\textsc{Step 6:} \textsl{If  $K_2(a)= \{1,2\}$, then $K_1(a)=\emptyset$
or both $n=3$ and $K_1(a)= \{2\}$.}

\textsc{Proof:} Let $K_2(a)= \{1,2\}$. Then $n\ge 3$. 
For $1\in K_1(a)$ we get 
\[
 g_n(a) &=& g_1^1(a)\,\,=\,\, a_1-\sum_{i\neq 1} a_i\,\,\le \,\,
 a_1-a_2 \,\,=\,\, -2g^2_2(a) \,\,=\,\, -2 g_n(a)
\]
contradicting Step 1. For $3\le k\in K_1(a)$ we get 
\[
 g_n(a)  \,\,=\,\, a_k -\sum_{i\neq k} a_i 
 &\le &a_k -a_1-a_2 \\
 &=& a_k -4g^2_1(a)-2g^2_2(a)
 \,\,=\,\, a_k- 6g_n(a)
\]
contradicting \eqref{lower bound gn} and $a_k>0$. If $K_1(a)= \{2\}$,
then 
\[
 g_n(a) &=&  g^2_2(a) \,\,=\,\, 
 g^1_2(a)    \,\,=\,\,2 g^2_2(a)-\sum_{i>2 }a_i \,\, =\,\,
 2 g_n(a) -  \sum_{i>2 }a_i 
\]
yields $a_k \le \sum_{i>2 }a_i = g_n(a)$ for $k\ge 3$,
and hence $K_3(a)=K_4(a)=\emptyset$, implying $n=3$ in view of 
\eqref{card.K1.4}.

\textsc{Step 7:} $K_4(a)= \{3,\ldots,n-1\}$.

\textsc{Proof:} This is trivial if $n\le 3$. 
For $n\ge 4$, inequality \eqref{card.K1.4} and Steps 2, 3, 5 and 6
yield $n\le 3+ \sharp K_4(a)$ and hence the claim.

\textsc{Step 8:} \textsl{$K_3(a)\subset\{n-1,n\}$.}

\textsc{Proof:} If $k\in\{3,\ldots,n-2\}$, then Step~7 implies
$k,k+1\in K_4(a)$ and thus 
$g^3_k(a)< \frac 13(a_k-a_{k+1}) = \frac 43(g^4_k(a) - g^4_{k+1}(a))=0$,
hence $k\notin  K_3(a)$.

\textsc{Step 9:} \textsl{If $n\ge 4$, then $K_1(a)=\emptyset$.}

\textsc{Proof:} Let $n\ge 4$ and $k\in K_1(a)$. Then $a_k = \max_{i=1}^n a_i$.
If $\ell\in\{3,\ldots,n-1\}$, then Step 7 yields $4g_n(a)=a_\ell$
and, using  \eqref{lower bound gn}, we get  $a_k\le a_\ell$
and thus $k=\ell$, for otherwise $g^1_k(a)\le a_k-a_\ell \le 0$.
For  $n\ge 5$ we thus get a contradiction by considering
$\ell=3$ and $\ell=4$. 

So let $n=4$.
Then $K_1(a)=\{3\}$, hence
$\sharp K_2(a) = 1$ by Steps 4 to 6, 
$K_3(a)\subset\{3,4\}$ by Step 8,
and $K_4(a)=\{3\}$ by Step 7.
If $4\in K_3(a)$, then $g^3_4(a)=g_n(a)=g^4_3(a)$ yields
$a_4=\frac 34 a_3$ and hence the contradiction
\[
 g_n(a) &=& g^1_3(a)\,\,< \,\,a_3-a_4\,\,=\,\,\frac{a_3}4\,\,=\,\,g_n(a)
\]
Thus, in view of \eqref{card.K1.4}, we must have $K_3(a) =\{3\}$. 
Hence $g^1_3(a)=g^3_3(a)$ and thus $a_1+a_2=\frac 23(a_3-a_4)$, 
so that  
for $\ell \in\{1,2\}$ we get $g^2_\ell(a) < \frac{a_1+a_2}2 = g^3_3(a)=g_n(a)$.
Hence $K_2(a)=\{3\}$. 
So $G:=  \{g^{\nu}_k \,:\, \nu\in\{1,2,3,4\},\,k\in K_\nu(a)\}
=\{g^1_3,g^2_3,g^3_3,g^4_3\}$, but the identity  
$g^1_3 -2g^2_3 -3 g^3_3 +4g^4_3=0$ shows that 
$G$ does not contain $n=4$ affinely independent functions,
contrary to our assumption preceding \eqref{card.K1.4}.

\textsc{Step 10:} \textsl{If $n\ge 4$, 
then \eqref{a_1/2...a_n/3} or \eqref{a_1/2...a_n}.}

\textsc{Proof:} The assumption together with \eqref{card.K1.4}
and Steps~3, 5, and 7 to 9 imply $K_1(a)=\emptyset$, $K_2(a)=\{1,2\}$, 
$K_4(a)=\{3,\ldots,n-1\}$, and $K_3(a)=\{n-1\}$ or $K_3(a)=\{n\}$.
The second possibility yields \eqref{a_1/2...a_n/3}, the first
yields \eqref{a_1/2...a_n}.

\textsc{Step 11:} 
\textsl{If   $n=3$, then \eqref{a_1/2...a_n/3} or \eqref{a_1/2...a_n}.}

\textsc{Proof:} Let $n=3$. Then $K_4(a)=\emptyset$, $K_3(a)\subset \{3\}$
and, using Step~14, $\emptyset \neq K_2(a) \subset \{1,2\}$. 
If $K_3(a)=\emptyset$, then \eqref{card.K1.4} and Steps 2, 5 and 6 yield
$K_1(a)=\{2\}$ and $K_2(A)=\{1,2\}$, and hence  \eqref{a_1/2...a_n}.
So let  $K_3(a)=\{3\}$. Then $K_1(a)=\emptyset$,
for otherwise either $3\in K_1(a)$ yielding 
$\frac{a_3}3 = g_n(a)= a_3-(a_1+a_2)$ and hence the contradiction
$g_n(a)= \frac{a_3}3 = \frac{a_1+a_2}2 > \max\{\frac {a_1}2,\frac{a_2-a_2}2\}
= g_n(a)$, or there is a $k\in \{1,2\}\cap K_1(a)$ yielding 
$\frac{a_3}3 = g_n(a)=|a_2-a_1| - a_3 $ and hence the contradiction
$g_n(a)= \frac{a_3}3 = \frac 14 |a_2-a_1| 
\le \frac 12 \max\{\frac {a_1}2,\frac{a_1+a_2}2\} =\frac 12 g_n(a)$.
Thus  \eqref{card.K1.4} yields $K_2(a)=\{1,2\}$  and \eqref{a_1/2...a_n/3}
follows.

\textsc{Step 12:} 
\textsl{If   $n=2$, then \eqref{a_1/2...a_n/3} or \eqref{a_1/2...a_n}.}

\textsc{Proof:} Let $n=2$. 
Then $K_3(a)=K_4(a)=\emptyset$ and Step~4 yields $K_2(a)=\{1\}$.
Thus \eqref{card.K1.4} and Step~2 yield either 
$K_1(a)=\{1\}$ and hence  \eqref{a_1/2...a_n},
or $K_1(a)=\{2\}$ and hence  \eqref{a_1/2...a_n/3}.
 
This completes our proof that  \eqref{a_1/2...a_n/3} or \eqref{a_1/2...a_n}  
holds. To prove $f_n(a)\le g_n(a)$ 
in each of these cases, we will consider  simple modifications of the 
martingales used in the proof of Theorem~\ref{Thm MG case}. 
Let $p$, $Y_1,\ldots,Y_n$, $S_0,\ldots,S_{n-1}$ and $X_1,\ldots,X_{n-1}$ 
be as in  \eqref{The martingale} and its preceding  four lines.

In case of \eqref{a_1/2...a_n/3} we define
\[
   S_n := \big(S_{n-1}Y_n \big)_+ &\text{ and }&
   X_n := S_n-S_{n-1} = (Y_n-1)(S_{n-1})_+ + (S_{n-1})_-
\]
Then $(S_i)_{i=0}^n$ is a submartingale and we have
$\E|S_n| =\E Y_n \E (S_{n-1})_+ =\frac 12$,  
$\E|X_1|=1$, $\E|X_i|= 2-2p$ for $i\in\{2,\ldots,n-1\}$,
and 
\[
 \E| X_n| 
  &=& \E|Y_n-1|\,\E(S_{n-1})_+ + \E(S_{n-1})_-  
 \,\,=\,\, (2-2p)\frac 12 +\frac 12 
 \,\,=\,\, \frac 32 -p
\]
Thus
\[
 f_n(1,2-2p,\ldots,2-2p,\frac 32 -p)&\le& \frac 12
\]
and hence $f_n(1,2,\ldots,2,\frac 32)\le \frac 12$
yielding $f_n(a) \le \frac {a_1}2 =g_n(a)$.

In case of \eqref{a_1/2...a_n} we define
\[
   S_n := \big(S_{n-1}\big)_+  &\text{ and }&
   X_n := S_n-S_{n-1} = \big(S_{n-1}\big)_-
\]
Then $(S_i)_{i=0}^n$ is a submartingale with 
$\E|S_n| = \frac 12 \E|S_{n-1}| =\frac 12$,  
$\E|X_1|=1$, $\E|X_i|= 2-2p$ for $i\in\{2,\ldots,n-1\}$,
and $\E |X_n| =\frac 12 \E|S_{n-1}| =\frac 12$. Thus
\[
 f_n(1,2-2p,\ldots,2-2p,\frac 12)&\le& \frac 12
\]
and hence $f_n(1,2,\ldots,2,\frac 12)\le \frac 12$
yielding again $f_n(a) \le \frac {a_1}2 =g_n(a)$.
\end{proof}

\end{document}